\documentclass[10pt]{amsart}
\usepackage[T1]{fontenc}
\usepackage{a4wide}
\usepackage[latin1]{inputenc}
\usepackage{amsmath}
\usepackage{amsfonts}
\usepackage{amssymb}
\usepackage{amsthm}
\usepackage{subfigure}
\newtheorem{theo}{Theorem}[section]
\newtheorem{prop}[theo]{Proposition}
\newtheorem{coro}[theo]{Corollary}
\newtheorem{lemm}[theo]{Lemma}

\theoremstyle{definition}

\theoremstyle{remark}
\newtheorem{rema}[theo]{Remark}
\newcommand{\Op}{\operatorname{Op}}

\title{Remarks on quantum ergodicity}
\author{Gabriel Rivi\`ere}
\date{\today}
\address{Laboratoire Paul Painlev\'e (U.M.R. CNRS 8524), U.F.R. de Math\'ematiques, Universit\'e Lille 1, 59655 Villeneuve d'Ascq Cedex, France}
\email{gabriel.riviere@math.univ-lille1.fr}

\begin{document}

\begin{abstract}
We prove a generalized version of the Quantum Ergodicity Theorem on smooth compact Riemannian manifolds without boundary. We apply it to prove some asymptotic properties on the distribution of typical eigenfunctions of the Laplacian in geometric situations where the Liouville measure is not (or not known to be) ergodic.
\end{abstract}

\maketitle

\section{Introduction}

Let $M$ be a smooth, compact, connected Riemannian manifold of dimension $d$ (without boundary). Denote by $L$ the normalized Liouville measure on the unit cotangent bundle $S^*M$ and by $g^t$ the geodesic flow on $S^*M$. Let $(\psi_j)_{j\in\mathbb{N}}$ be an orthonormal basis of eigenfunctions of $-\Delta_g$ associated to a nondecreasing sequence of eigenvalues $(\lambda_j^2)_{j\in\mathbb{N }}$, i.e.
$$-\Delta_g\psi_j=\lambda_j^2\psi_j,\ \|\psi_j\|_{L^2(M)}=1.$$

In the following, we will write $N(\lambda):=\sharp\{j:\lambda_j^2\leq\lambda^2\}.$ Our goal in this note is to describe the asymptotic distribution of this sequence of eigenfunctions as $\lambda_j^2$ tends to infinity. For that purpose, we introduce the following distribution on $S^*M$:
$$\forall a\in\mathcal{C}^{\infty}(S^*M),\ \mu_{j}(a)=\int_{S^*M}ad\mu_{j}:=\langle\psi_j,\Op(a)\psi_j\rangle,$$
where $\Op(a)$ is a pseudodifferential operator with principal symbol $a$. It is a classical fact to check that any accumulation point\footnote{Convergence is for the standard topology on $\mathcal{D}'(S^*M)$.} of this sequence belongs to the set $\mathcal{M}(S^*M,g^t)$ of $(g^t)_t$-invariant probability measures on $S^*M$~\cite{Ge91, Bu97, Zw12}.

If the Liouville measure is ergodic, the Quantum Ergodicity Theorem states that there exists a subset $S$ of density\footnote{Recall that $S\subset \mathbb{N}$ have density $1$ if $\displaystyle\lim_{n\rightarrow+\infty}\frac{1}{n}\sharp\{j\in S:1\leq j\leq n\}=1.$} $1$ in $\mathbb{N}$ such that the sequence $(\mu_j)_{j\in S}$ converges to the Liouville measure $L$~\cite{Sh74, Ze87, CdV85}. In other words, it means that the eigenfunctions become equidistributed in $S^*M$. We refer the reader to~\cite{Ze10} for a recent detailed survey on related issues. In our context, the main example of application is given by geodesic flows on manifolds of negative curvature or more generally by uniformly hyperbolic geodesic flows: in this setting, the Liouville measure is known to be ergodic.

Here, we are interested in the case where we drop the ergodicity assumption. The main examples we have in mind are geodesic flows for which the Liouville measure is ergodic on a subset of positive measure. For instance, this kind of situations occurs when the geodesic flow is supposed to be nonuniformly hyperbolic~\cite{BaPe06}. In these cases, we derive properties on the asymptotic distributions of the eigenmodes. For that purpose, we prove an alternative version of the Quantum Ergodicity Theorem that does not rely on ergodicity. Then, we apply this result in several geometric contexts. For example, we obtain a kind of equidistribution property for subsequences of eigenfunctions of the Laplacian on surfaces of nonpositive curvature with genus $\geq 2$. In this setting, the Liouville measure is not known to be ergodic; thus, the standard Quantum Ergodicity Theorem does not apply a priori.

\section{Statement of the main result}
\label{s:statement}

Let $\Lambda$ be a subset of $S^*M$. We say that it satisfies the Birkhoff property if
\begin{equation}\label{e:birkhoff}\forall\rho\in\Lambda,\ \exists\ L_{\rho}\in\mathcal{M}(S^*M,g^t),\ \lim_{T\rightarrow+\infty}\frac{1}{T}\int_0^T\delta_{g^t\rho}dt =L_{\rho},\end{equation}
where the convergence is for the weak-$\star$ topology, i.e.
$$\forall a\in \mathcal{C}^0(S^*M),\ \forall\rho\in \Lambda,\ \lim_{T\rightarrow+\infty}\frac{1}{T}\int_0^Ta(g^t\rho)dt =L_{\rho}(a).$$
For every subset $\Lambda\subset S^*M$ satisfying property~\eqref{e:birkhoff}, we introduce $\text{Cv}(\Lambda)$ which is the closure in $\mathcal{D}'(S^*M)$ of the convex hull of $\{L_{\rho}:\rho\in \Lambda\}$. We emphasize that the set $\text{Cv}(\Lambda)$ depends on the choice of $\Lambda$ and that any element in $\text{Cv}(\Lambda)$ is a probability measure on $S^*M$ invariant under the geodesic flow.

Thanks to the Birkhoff Ergodic Theorem~\cite{EiWa11} -section $6.1$, there exists $\Lambda\subset S^*M$ satisfying~\eqref{e:birkhoff} and $L(\Lambda)=1$. The main result of this note is the following version of the Quantum Ergodicity Theorem:

\begin{theo}\label{t:GQE} Let $\Lambda\subset S^*M$ satisfying
property~\eqref{e:birkhoff} and such that $L(\Lambda)=1$. Let $(\psi_j)_{j\in\mathbb{N}}$ be an orthonormal basis of eigenfunctions of $-\Delta_g$ associated to a nondecreasing sequence of eigenvalues $(\lambda_j^2)_{j\in\mathbb{N }}$, i.e.
$$-\Delta_g\psi_j=\lambda_j^2\psi_j,\ \|\psi_j\|_{L^2(M)}=1.$$

Then, there exists $S\subset\mathbb{N}$ of density $1$ such that any accumulation point of the sequence $(\mu_j)_{j\in S}$ belongs to $\operatorname{Cv}(\Lambda)$.
\end{theo}

\begin{rema}

We emphasize that the result holds for any choice of $\Lambda\subset S^*M$ of full measure satisfying the Birkhoff property~\eqref{e:birkhoff} (the result is of course more interesting when the set $\text{Cv}(\Lambda)$ is not too big). We also underline that this theorem is true for any orthonormal basis of eigenfunctions of $\Delta_g$ and that we do not need to make any particular assumption on the manifold (like ergodicity for instance). Even if this generalization is quite natural, we did not find any trace of such a result in the literature. 

\end{rema}

If $M$ is the sphere $\mathbb{S}^d$ endowed with its canonical metric, the set $\text{Cv}(\Lambda)$ is always equal to $\mathcal{M}(S^*M,g^t)$ for $\Lambda$ of full measure. Thus, the result is empty as we already know that any accumulation point of the sequence $(\mu_j)_{j\geq 0}$ is an invariant probability measure. In the ``opposite'' case where the Liouville measure is ergodic for the geodesic flow, we recover the standard Quantum Ergodicity Theorem~\cite{Sh74, Ze87, CdV85} as, thanks to the Birkhoff Ergodic Theorem, one can pick $\Lambda$ of full measure satisfying~\eqref{e:birkhoff} and $L_{\rho}=L$ on $\Lambda$ (and thus $\text{Cv}(\Lambda)=\{L\}$). In section~\ref{s:examples}, we will provide other examples of applications.\\

In the physics literature, the ``semiclassical eigenfunctions hypothesis'' states that the eigenmodes $(\psi_j)_{j\geq 0}$ must be asymptotically concentrated into regions of phase space which a typical orbit explores in the long time limit~\cite{Per73, Be77}. In our context, the set $\Lambda$ could represent in some sense a set of typical orbits and the measure $L_{\rho}$ is the canonical measure associated to the orbit of a point $\rho$ in the phase space $S^*M$. 

Regarding this conjecture, it seems natural to understand when there exist a subset $\Lambda$ of full measure and a typical subsequence of eigenmodes $(\psi_j)_{j\in S}$ such that the accumulation points of $(\mu_j)_{j\in S}$ are exactly given by the closure of $\{L_{\rho}:\rho\in \Lambda\}$. This question was for instance raised by Shnirelman in~\cite{Sh93} -- end of paragraph $AD.2$. In such generality, it is not true as there exist geometric situations where the set $\{L_{\rho}:\rho\in \Lambda\}$ cannot be reduced to $L$ while there exists a typical family of states that converges to $L$~\cite{Ze92, Ze10, Gu09, Ze12}.

Our theorem shows that, for a typical family $(\mu_j)_{j\in S}$, the accumulation points belong to a larger set than the closure of $\{L_{\rho}:\rho\in \Lambda\}$, precisely they belong to the closure of its convex hull. We underline that results related to these questions were also obtained by Marklof and O'Keefe for specific families of quantum maps with divided phase space~\cite{MaOk05}, and by Galkowski for Hamiltonian flows with divided phase space on manifolds with piecewise smooth boundary~\cite{Ga12}.

We will explain in paragraph~\ref{s:theo} how one can get our generalized version of the Quantum Ergodicity Theorem by implementing an idea used by Sj\"ostrand in the context of damped wave equations~\cite{Sj00}. In fact, our proof will combine Hahn-Banach Theorem with the following main lemma that makes the connection with the results in~\cite{Sj00} more explicit:


\begin{lemm}\label{l:GQE} Let $a$ be an element in $\mathcal{C}^{\infty}(S^*M,\mathbb{R})$. Then, there exists $S\subset\mathbb{N}$ of density $1$ such that,
$$\operatorname{essinf}\ L_{\rho}(a)\leq\liminf_{j\rightarrow+\infty,j\in S}\mu_j(a)\leq\limsup_{j\rightarrow+\infty,j\in S}\mu_j(a)\leq\operatorname{esssup}\ L_{\rho}(a).$$

\end{lemm}

As in the classical proof of the Quantum Ergodicity Theorem, the main ingredient of the proof of this lemma is a local Weyl law combined to an averaging of $a$ along geodesics and to a positivity argument. 

\begin{rema}\label{r:esssup} By construction, $\rho\mapsto L_{\rho}(a)$ belongs to $L^{\infty}(S^*M,dL)$, and thus$|\operatorname{essinf}L_{\rho}(a)|,\ |\operatorname{esssup}L_{\rho}(a)|<+\infty$. Compared with our main theorem, these two quantities do not depend on any choice of subset $\Lambda$. For any subset $\Lambda$ of full measure satisfying~\eqref{e:birkhoff}, one has
$$\inf_{\rho\in\Lambda}L_{\rho}(a)\leq\operatorname{essinf}L_{\rho}(a)\leq\operatorname{esssup}L_{\rho}(a)\leq\sup_{\rho\in\Lambda}L_{\rho}(a).$$  
\end{rema}

\subsection*{Organization of the following} 
In the next paragraph, we apply our theorem in several geometric situations. After that, we give the proof of theorem~\ref{t:GQE} and provide the proof of an intermediary proposition that we used in our applications.

In the following, we will sometimes use the notation $\mathcal{D}(S^*M)$ for the space $\mathcal{C}^{\infty}(S^*M)$ when we want to emphasize that we are working with distributions.

\section{Application of theorem~\ref{t:GQE}}
\label{s:examples}
Before entering the details of the proof, we describe several geometric situations where the Liouville measure is a priori not supposed to be ergodic and we apply theorem~\ref{t:GQE} in order to derive some weak equidistribution properties for the eigenfunctions.

\subsection{Geodesic flows with divided phase space}
\label{ss:Donnay}

A direct consequence of theorem~\ref{t:GQE} is the following property:
\begin{coro}\label{c:maincoro} Suppose there exist $\Lambda'$ of positive Liouville measure, $I$ finite and a family $(\Lambda_i)_{i\in I}$ of invariant subsets such that 
\begin{itemize}
\item $L(\Lambda_i)>0$ for every $i\in I$;
\item $\Lambda_i\cap\Lambda_j$ is empty when $i\neq j$;
\item $\cup_i\Lambda_i=\Lambda'$;
\item $L_{|\Lambda_i}$ is ergodic for every $i\in I$.
\end{itemize}
Then, for every $\Lambda_0\subset\Lambda'^c$ satisfying~\eqref{e:birkhoff} and $L(\Lambda_0)=L(\Lambda'^c)$, there exists a subset $S\subset\mathbb{N}$ of density $1$ such that any accumulation point of the sequence $(\mu_j)_{j\in S}$
is of the form
\begin{equation}\label{e:donnay}
\mu=\alpha\sum_{i\in I}t_i\frac{L_{|\Lambda_i}}{L(\Lambda_i)}+(1-\alpha)\nu_{0},
\end{equation}
where $0\leq \alpha\leq 1$, $\displaystyle 0\leq t_i\leq 1$, $\sum_{i\in I}t_i=1$ and $\nu_{0}$ belongs to $\operatorname{Cv}(\Lambda_0)$. 

\end{coro}
 
In the case where $L(\Lambda')=1$, this corollary tells us that, along a subsequence $S$ of density $1$, any accumulation point of the sequence $(\mu_j)_{j\in S}$ must be absolutely continuous with respect to $L$. Moreover, it gives us an explicit expression for the accumulation points.
Galkowski recently proved a similar result in the case of manifolds with piecewise smooth boundaries under the assumption that $\Lambda'$ satisfies $L_{|\Lambda'}$ ergodic, $L(\Lambda')>0$ and $L(\partial \Lambda'-\Lambda')=0$~\cite{Ga12}. It would be interesting to understand if theorem~\ref{t:GQE} could also be obtained in the setting of manifolds with piecewise smooth boundaries.



\begin{rema}\label{r:pesin} Even if there is no hypothesis on nonuniform hyperbolicity in the above statement, we underline that the assumptions of the corollary appear naturally in this context -- we refer the reader to~\cite{BaPe06} for recent lecture notes on this topic. Using Pesin classical work on nonuniform hyperbolicity~\cite{Pe77} (see also Theorem $11.5$ in~\cite{BaPe06}), we recall that, if there exists an invariant subset $\Lambda\subset S^*M$ of positive Liouville measure such that \emph{the geodesic flow is nonuniformly hyperbolic on} $\Lambda$, then there exists $I$ \emph{at most countable} and a family $(\Lambda_i)_{i\in I}$ of invariant subsets such that 
\begin{itemize}
\item $\cup_{i\in I}\Lambda_i=\Lambda$ (modulo a subset of $0$ Liouville measure);
\item $L(\Lambda_i)>0$ for every $i\in I$;
\item $\Lambda_i\cap\Lambda_j$ is empty when $i\neq j$;
\item $L_{|\Lambda_i}$ is ergodic for every $i\in I$.
\end{itemize}
\end{rema}

\begin{rema}
In~\cite{Gu09}, Gutkin constructs a nonuniformly hyperbolic geodesic flow on a billiard table for which the Liouville measure is not ergodic. However, his system presents 
a symmetry that allows him to obtain a stronger result than corollary~\ref{c:maincoro}, precisely he proves the existence of a subsequence of density $1$ converging to 
the Liouville measure. It is not clear to the author whether there exist smooth compact Riemannian manifolds without boundary which satisfy the assumption of the corollary 
with  $L(\Lambda')=1$ and $|I|\geq 2$, and which do not present a symmetry like the one in~\cite{Gu09}. 
We also do not know if ergodicity or even $I$ finite are ``generic properties'' in the family of smooth nonuniformly hyperbolic geodesic flows.  
\end{rema}

\begin{rema} In order to provide a concrete example, we briefly discuss a construction due to Donnay~\cite{Do88} -- section $11$. In this reference, he constructs Riemmanian metrics on the sphere $\mathbb{S}^2$ for which the phase space splits into a chaotic component and a regular 
one. His strategy is to remove three or more points from the sphere and to endow the induced punctured surface with the Poincar\'e metric; then, he attaches 
smoothly a so called ``light-bulb cap'' in a neighborhood of each deleted point. \emph{The geodesic flow he obtains is not ergodic for $L$}. The phase space $S^*M$ contains 
two disjoint invariant subsets $\Lambda_{\text{chaotic}}$ and $\Lambda_{\text{regular}}$ of \emph{positive Liouville measure} satisfying 
$L(\Lambda_{\text{chaotic}}\cup\Lambda_{\text{regular}})=1$. More precisely, $\Lambda_{\text{regular}}$ consists of orbits that stay in the caps 
and $g^t_{|\Lambda_{\text{chaotic}}}$ is nonuniformly hyperbolic. Thanks to remark~\ref{r:pesin}, the corollary applies here.

\end{rema}





\subsection{More on divided phase space}

The statements of the previous paragraph do not forbid that the ergodic component associated to a subset $\Lambda_i$ of positive Liouville measure has a weight $t_i=0$. This motivates the following proposition which gives partial informations on this question.

\begin{prop}\label{p:genprop} Suppose there exists an invariant subset $\Lambda'$ such that $\Lambda'$ contains a nonempty open ball (modulo a subset of $0$ Liouville measure) and such that $L_{|\Lambda'}$ is ergodic.\\
Then, there exists $\Lambda_1\subset\Lambda'^c$ satisfying $L(\Lambda_1)=L(\Lambda'^c)$ such that,
for every $\Lambda_0\subset\Lambda_1$ verifying~\eqref{e:birkhoff} and $L(\Lambda_1)=L(\Lambda_0)$ and for every $0\leq \delta<L(\Lambda')$, there exists $S_{\delta}\subset \mathbb{N}$ of density\footnote{By density $\geq\delta$, we mean that $\displaystyle\liminf_{\lambda\rightarrow+\infty}\frac{1}{N(\lambda)}\sharp\{j:\lambda_j^2\leq\lambda^2\ \text{and}\ j\in S_{\delta}\}\geq\delta.$}  $\geq\delta$ such that any accumulation point of the sequence $(\mu_j)_{j\in S_{\delta}}$ is of the form~\eqref{e:donnay} with
$$0<\frac{L(\Lambda')-\delta}{1-\delta}\leq\alpha\leq 1.$$
\end{prop}

Even if the conclusion of the proposition is slightly technical, we emphasize that we only need to assume that $\Lambda'$ contains a nonempty open ball (modulo a subset of $0$ Liouville measure) in order to have positive mass on the ergodic component $\Lambda'$.

This result is true for any orthonormal basis of eigenfunctions of $\Delta_g$. It can be obtained as an application of lemma~\ref{l:GQE} -- see paragraph~\ref{s:prop} for details. It provides a sufficient condition to observe a kind of equidistribution property on an invariant subset $\Lambda'$ -- see~\cite{Sh93} for related questions.


\begin{rema}\label{r:generalization} The hypothesis is slightly different from the one of corollary~\ref{c:maincoro} but the argument could be adapted to treat the case where the assumptions of the corollary are satisfied with $\Lambda'$ containing a nonempty open ball (modulo a subset of $0$ Liouville measure) -- see the end of paragraph~\ref{s:prop} for details.
Precisely, we will check that there exists $\Lambda_1\subset\Lambda'^c$ satisfying $L(\Lambda_1)=L(\Lambda'^c)$ such that,
for every $\Lambda_0\subset\Lambda_1$ verifying~\eqref{e:birkhoff} and $L(\Lambda_1)=L(\Lambda_0)$, there exists a subset $S'$ of positive density such that 
any accumulation point of the sequence $(\mu_j)_{j\in S'}$ is of the form given by equation~\eqref{e:donnay} with $\alpha>0$. 
\end{rema}

\subsection{Surfaces of nonpositive curvature}

We will now give an application of the previous proposition in the context of nonpositively curved manifolds. We suppose that $M$ is a \emph{surface of nonpositive curvature} of genus $\geq 2$. For any point $x$ in $M$, we will denote by $K(x)\leq 0$ the sectional curvature at point $x$. We will make a small abuse of notations and use also the notation $K$ for its canonical lift on $S^*M$. Following~\cite{BaPe06}, we introduce the following \emph{invariant} subset
\begin{equation}\label{e:pesin}
\Lambda':=\left\{\rho\in S^*M:\limsup_{T\rightarrow+\infty}\frac{1}{T}\int_0^TK\circ g^t(\rho)dt<0\right\}.
\end{equation}                                                                                                            
According to~\cite{BaPe06} (Theorems $2.4$ and $17.7$), the set $\Lambda'$ is \emph{open (modulo a subset of $0$ Liouville measure)} and it satisfies $L(\Lambda')>0$ and $L_{|\Lambda'}$ is ergodic. In particular, the subset $\Lambda'$ satisfies all the requirements of proposition~\ref{p:genprop}. It is also shown in this reference that $\Lambda'$ \emph{is everywhere dense} but we will not use this fact here. 

We underline that it is still an open question to determine whether $L(\Lambda')=1$ or not for any surface of nonpositive curvature of genus $\geq 2$. In other words, \emph{it is not known if the Liouville measure is ergodic or not}. Yet, we can prove that the eigenfunctions satisfy some equidistribution properties in this negatively curved part of the surface. For that purpose, we observe that, thanks to Birkhoff Ergodic Theorem, $L_{\rho}(K)$ is well defined a.e. on $\Lambda'^c$. By definition of $\Lambda'$, we obtain 
$$L_{\rho}(K)=0\ \text{a.e. on}\ \Lambda'^c.$$ 

As the subset $\Lambda'$ satisfies all the requirements of proposition~\ref{p:genprop}, the following corollary holds:


\begin{coro}\label{c:mainresult} If $M$ is a surface of nonpositive curvature $K(x)$ and of genus $\geq 2$, then there exists $\Lambda_1\subset \{L_{\rho}(K)=0\}\cap\Lambda'^c$ (with $L(\Lambda_1)=L(\Lambda'^c)$) such that, for any $\Lambda_0\subset\Lambda_1$ satisfying~\eqref{e:birkhoff} and $L(\Lambda_0)=L(\Lambda'^c)$, and, for every $0\leq \delta<L(\Lambda')$,
\begin{enumerate}

\item there exists a subset $S$ of density $1$ in $\mathbb{N}$ such that any accumulation point of the sequence $(\mu_j)_{j\in S}$ is of the form
$$\alpha \frac{L_{|\Lambda'}}{L(\Lambda')}+(1-\alpha)\nu_0,$$
where $0\leq\alpha\leq 1$ and $\nu_0$ belongs to $\operatorname{Cv}(\Lambda_0)$;
\item there exists $S_{\delta}\subset S$ of density $\geq\delta$ such that, for any accumulation point of the sequence $(\mu_j)_{j\in S_{\delta}}$, one has 
$$0<\frac{L(\Lambda')-\delta}{1-\delta}\leq\alpha\leq 1.$$
\end{enumerate}


\end{coro}

This result is true for any orthonormal basis of eigenfunctions of $\Delta_g$. This corollary tells us that a positive proportion of eigenmodes are asymptotically equidistributed in $\Lambda'$ even if we do not have ergodicity of the Liouville measure on the entire phase space. We emphasize that the assumptions on $\Lambda_0$ implies that $\nu_0(K)=0$ for any $\nu_0$ in $\text{Cv}(\Lambda_0)$.




\begin{rema}
If we project the distributions $\mu_j$ on the base, we find the following notable consequence of the previous corollary: any accumulation point of the sequence $(K|\psi_j|^2\text{vol}_M)_{j\in S}$ is of the form $cK\text{vol}_M$ where $c\geq 0$ is a constant. Moreover, thanks to part~$(2)$ of the corollary, there exist subsequences of positive density for which $c$ can be chosen positive.
\end{rema}

\begin{rema} In the case where $\dim M\geq 2$, one can introduce the following subset of $S^*M$:
$$\Lambda':=\left\{\rho\in S^*M:\limsup_{T\rightarrow+\infty}\frac{1}{T}\int_0^TK_{\pi\circ g^t\rho}(g^t\rho,g^t\rho')dt<0,\ \text{for every}\ \rho'\ \text{orthogonal to}\ \rho\right\},$$
where $\pi:S^*M\rightarrow M$ is the canonical projection on $M$ and $K_x(v_1,v_2)$ is the sectional curvature for $x$ in $M$ and $v_1,v_2$ in $T_x^*M$. Suppose now that $M$ has nonpositive curvature, i.e. $K_x(v_1,v_2)\leq 0$ for every $x$ in $M$ and every $v_1,v_2$ in $T_x^*M$. Under some extra geometric assumptions\footnote{For more details on these assumptions, we refer the reader to~\cite{BaPe06}, sections $2.7$ and $17.1$.} on $M$ that are always satisfied by nonpositively curved surfaces of genus $\geq 2$, the set $\Lambda'$ is again open (modulo $0$) and everywhere dense and it satisfies $L(\Lambda')>0$ and $L_{|\Lambda'}$ is ergodic. Then, corollary~\ref{c:mainresult} can be extended in $\dim M\geq 2$ modulo the above extra geometric assumptions. 
\end{rema}

\subsection{Flat torus} In this paragraph, we apply theorem~\ref{t:GQE} to the flat torus $\mathbb{T}^{d}$ for which there is no subset $\Lambda'$ of positive measure on which $L$ is ergodic. We introduce the subset of ``irrational'' vectors 
$$\Lambda:=\left\{(x,\xi)\in S^*\mathbb{T}^d:\ \forall p\in\mathbb{Z}^d-\{0\},\ p.\xi\neq 0\right\}.$$
This set has full Liouville measure and $\Lambda$ satisfies property~\eqref{e:birkhoff} with $L_{\rho}=dx\times\delta_{\xi}$. In particular, the projection on $\mathbb{T}^d$ of any element in $\text{Cv}(\Lambda)$ is the Lebesgue measure $dx$. Applying theorem~\ref{t:GQE}, we obtain the following corollary:
\begin{coro}\label{c:torus} For any orthonormal basis $(\psi_{j})_{j\in\mathbb{N}}$ of eigenfunctions of $\Delta$ on $\mathbb{T}^d$, there exists a subset $S$ of density $1$ in $\mathbb{N}$ such that
$$\forall a\in\mathcal{C}^0(\mathbb{T}^d),\ \lim_{j\rightarrow+\infty, j\in S}\int_{\mathbb{T}^d}a(x)|\psi_j(x)|^2dx=\int_{\mathbb{T}^d}a(x)dx.$$

\end{coro}

This result is the analogue on $\mathbb{T}^d$ of Marklof-Rudnick's recent result on equidistribution of eigenfunctions on rational polygons~\cite{MaRu12}.

\section{Proof of the main result}\label{s:theo}

The proof follows classical ideas taken from~\cite{Sh74, Ze87, CdV85, Sj00} that we carefully combine to get our main theorem.

\subsection{Preliminary remarks}

Before getting into the details of the proof, we mention a few facts on the properties of the quantization procedure $\Op$. First, recall that it is not defined in a 
canonical way: it depends on a choice of coordinate charts and on a choice of quantization procedure on $\mathbb{R}^{2d}$ -- see section $II.5$ in~\cite{Ta} or~chapter~$14$ in~\cite{Zw12}. 
A standard choice is to take the Weyl quantization on $\mathbb{R}^{2d}$. Then, for two different choices of coordinate charts, we obtain two quantization procedures $\Op$ and $\Op'$. Yet, one can show that, for any $a$ in $\mathcal{C}^{\infty}(S^*M)$, $\Op(a)-\Op'(a)$ is a pseudodifferential operator of order $-1$. In particular, one has
$$\lim_{j\rightarrow+\infty}\langle\psi_j,(\Op(a)-\Op'(a))\psi_j\rangle_{L^2(M)}=0,$$
and thus the accumulation points of the sequence $(\mu_j)_{j\geq 0}$ do not depend on the choice of coordinates on $M$. Let us now briefly recall a few properties of this quantization procedure that we will use at different steps of our argument.\\

For every $a$ in $\mathcal{C}^{\infty}(S^*M)$, $\Op(a)$ defines a bounded operator in $L^2(M)$. Moreover, if $a$ is real valued, then $\Op(a)$ is selfadjoint. In particular, the sequence of distributions $\mu_j$ is real valued, i.e. $\mu_j(a)$ belongs to $\mathbb{R}$ when $a$ is real valued.

We introduce the average of $a\in\mathcal{C}^{\infty}(S^*M)$ at time $T>0$,
$$ a_T(\rho):=\frac{1}{T}\int_0^Ta\circ g^t(\rho)dt.$$
As $(\psi_j)_{j\geq 0}$ is a sequence of eigenmodes of $\Delta_g$, one can apply the Egorov Theorem (theorem~$15.2$ in~\cite{Zw12}) and get
\begin{equation}\label{e:egorov}
\forall T>0,\ \mu_j(a)=\mu_j(a_T)+o_T(1),
\end{equation}
where $o_T(1)$ is a remainder that depends on $a$ and $T$ and that tends to $0$ as $j$ tends to infinity.

We will also use the following local Weyl law (theorem~$15.3$ in~\cite{Zw12}):
\begin{equation}\label{e:weyl}
\forall b\in\mathcal{C}^{\infty}(S^*M),\ \frac{1}{N(\lambda)}\sum_{j:\lambda_j^2\leq\lambda^2}\mu_j(b)=\int_{S^*M}bdL+o(1),\ \text{as}\ \lambda\rightarrow+\infty,
\end{equation}
where the remainder depends on $b$.

An inconvenience of choosing the Weyl quantization on $\mathbb{R}^{2d}$ is that it does not satisfy a positivity property. This can be solved by using a different 
quantization on $\mathbb{R}^{2d}$. For instance, as in the proof of~\cite{CdV85} (section $1.1$), one can take the so-called Friedrichs quantization 
on $\mathbb{R}^{2d}$ -- see section~$VII.2$ in~\cite{Ta} for details. 
After choosing a family of coordinate charts, it gives us a quantization procedure $\Op^+$ on $M$ such that, for every $b\in\mathcal{C}^{\infty}(S^*M)$,
$$b\geq 0\Longrightarrow \Op^+(b)\geq 0.$$ 
The Weyl and the Friedrichs quantization are ``equivalent'' on $\mathbb{R}^{2d}$ -- Theorem~$2.2$ of Ch.~$VII$~\cite{Ta}. In particular, one has, for every $b\in\mathcal{C}^{\infty}(S^*M)$,
\begin{equation}\label{e:fried}
\lim_{j\rightarrow+\infty}\langle\psi_j,(\Op(b)-\Op^+(b))\psi_j\rangle_{L^2(M)}=0.
\end{equation}
Thus, the accumulation points of the sequence $\mu_j$ do not change if we replace $\Op$ by $\Op^+$ in the definition of $\mu_j$. Moreover, the invariance relation~\eqref{e:egorov} and the local Weyl law~\eqref{e:weyl} are still valid if we replace $\Op$ by $\Op^+$ in the definition of $\mu_j$.

\subsection{Proof of lemma~\ref{l:GQE}}

We start our proof by giving the proof of the main lemma~\ref{l:GQE}. Let $a$ be an element in $\mathcal{C}^{\infty}(S^*M,\mathbb{R})$. 
In order to simplify the presentation, denote $A_0:=\operatorname{esssup} L_{\rho}(a)$ -- see remark~\ref{r:esssup}. By definition, one has that, for every $\delta>0$,
\begin{equation}\label{e:esssup}L\left(\left\{\rho\in S^*M: a_T(\rho)\geq A_0+\delta\right\}\right)\rightarrow 0,\ \text{as}\ T\rightarrow+\infty.\end{equation}


Fix now $T>0$ and $\epsilon>0$. Properties~\eqref{e:egorov} and~\eqref{e:fried} tell us that, as $j$ tends to $\infty$,
$$\mu_j(a)=\mu_j(a_T)+o_T(1)=\langle\psi_j,\Op^+( a_T)\psi_j\rangle+o_T(1),$$
where the remainder depends on $T$. As in~\cite{Sj00}, one can define a new smooth function $\tilde{a}_T\leq a_T$ on $S^*M$ such that
\begin{itemize}
 \item $\tilde{a}_T(\rho)= a_T(\rho)$ when $ a_T(\rho)\leq A_0+\frac{\sqrt{\epsilon}}{2}$.
 \item $\tilde{a}_T(\rho)\leq A_0+\sqrt{\epsilon}$ otherwise.
\end{itemize}

\begin{rema}
In order to construct the function $\tilde{a}_T$, one can fix an increasing smooth function $\theta_{\epsilon}:\mathbb{R}\rightarrow\mathbb{R}$ which satisties
$$\forall v\in\mathbb{R},\ \theta_{\epsilon}(v)\leq v,\ |\theta'(v)|\leq 1,$$
and
$$\forall v\leq A_0+\frac{\sqrt{\epsilon}}{2},\ \theta_{\epsilon}(v)=v,\ \text{and}\ \theta_{\epsilon}(v)\leq A_0+\sqrt{\epsilon}\ \text{otherwise}.$$
Then, one defines $\tilde{a}_T=\theta_{\epsilon}(a_T)$. In particular, $\|\tilde{a}_T - a_T\|_{\infty}$ is bounded independently of $T>0$. 
\end{rema}
As $j$ tends to infinity, one has the following equality:
$$\mu_j(a)=\langle\psi_j,\Op^+(a_T-\tilde{a}_T)\psi_j\rangle+\langle\psi_j,\Op^+( \tilde{a}_T)\psi_j\rangle+o_T(1).$$

The idea of introducing this new function is taken from~\cite{Sj00} where it was used to study spectral asymptotics of the damped wave equation. In the following lines, we will show that
\begin{itemize}
\item most of the terms in the sequence $(\langle\psi_j,\Op^+(a_T-\tilde{a}_T)\psi_j\rangle)_{j\geq 0}$ are small following arguments from~\cite{Ze87, CdV85}; 
\item the other term in the RHS will be less than $A_0+\sqrt{\epsilon}+o_T(1)$ by construction.
\end{itemize}
A careful combination of these two facts will finally allow us to get our conclusion.\\

From~\eqref{e:fried} and the local Weyl law~\eqref{e:weyl}, one has
$$\frac{1}{N(\lambda)}\sum_{\lambda_j^2\leq\lambda^2}\langle\psi_j,\Op^+(a_T-\tilde{a}_T)\psi_j\rangle=\int_{S^*M}(a_T-\tilde{a}_T)dL+o_T(1),$$
where each term in the sum is nonnegative (as $a_T-\tilde{a}_T\geq 0$). Thanks to our construction, one has
$$\int_{S^*M}(a_T-\tilde{a}_T)dL\leq C_{a,\epsilon}L\left(\left\{\rho\in S^*M: a_T(\rho)\geq A_0+\frac{\sqrt{\epsilon}}{2}\right\}\right).$$
Fix $\eta>0$. Combining~\eqref{e:esssup} to the two previous relations, there exists $T_{\epsilon,\eta}>0$ such that, for every $T\geq T_{\epsilon,\eta}$, one can find $\lambda_T>0$ satisfying
$$\lambda\geq\lambda_T\Longrightarrow \frac{1}{N(\lambda)}\sum_{\lambda_j^2\leq\lambda^2}\langle\psi_j,\Op^+(a_T-\tilde{a}_T)\psi_j\rangle\leq \eta \epsilon.$$
We now fix $T=T_{\epsilon,\eta}.$ Denote $D_{\epsilon}(\lambda):=\{j:\lambda_j^2\leq\lambda^2\ \text{and}\ \langle\psi_j,\Op^+(a_T-\tilde{a}_T)\psi_j\rangle\geq\sqrt{\epsilon}\}$. 
Thanks to the Tchebychev inequality, we obtain that 
$\frac{\sharp D_{\epsilon}(\lambda)}{N(\lambda)}\leq\eta\sqrt{\epsilon}$ for $\lambda\geq\lambda_T$. This means that most of the terms in the sequence of nonnegative numbers $(\langle\psi_j,\Op^+(a_T-\tilde{a}_T)\psi_j\rangle)_{j\geq 0}$ are small.

As $\Op^+$ is positive and $\tilde{a}_T\leq A_0+\sqrt{\epsilon}$, one has that, for $\lambda_j^2$ larger than some $A>0$, the term $\langle\psi_j,\Op^+( \tilde{a}_T)\psi_j\rangle+o_T(1)$ is less than $A_0+2\sqrt{\epsilon}$. Thanks to the above discusion, we can write
$$\sharp\left\{j:\lambda_j^2\leq\lambda^2\ \text{and}\ \mu_j(a)\leq A_0+3\sqrt{\epsilon}\right\}\geq\sharp\left\{j:\lambda_j^2\leq A\ \text{and}\ \mu_j(a)\leq A_0+3\sqrt{\epsilon}\right\}$$
$$\hspace{5cm}+\sharp\left\{j:A\leq\lambda_j^2\leq\lambda^2\ \text{and}\ \langle\psi_j,\Op^+(a_T-\tilde{a}_T)\psi_j\rangle<\sqrt{\epsilon}\right\}.$$
If we denote $S_{\epsilon}:=\left\{j:\ \mu_j(a)\leq A_0+3\sqrt{\epsilon}\right\}$, then we have 
$$\liminf_{\lambda \rightarrow+\infty}\frac{\sharp \left(S_{\epsilon}\cap\left\{j:\lambda_j^2\leq\lambda^2\right\}\right)}{N(\lambda)}\geq 1-\eta\sqrt{\epsilon}.$$ 
This is true for any $\eta>0$ which implies that $S_{\epsilon}$ has density $1$.

Using the procedure of paragraph~$5$ in~\cite{CdV85} -- see remark~\ref{r:cdv} below, one can then obtain a subset $S_{0}\subset\mathbb{N}$ of density $1$, such that any accumulation point of the sequence $(\mu_j(a))_{j\in S_0}$ is~$\leq A_0$. This achieves the proof of the upper bound in lemma~\ref{l:GQE} and the lower bound can be easily derived by considering $-a$.

\begin{rema}\label{r:cdv}
 For the sake of completeness, we reproduce the argument used in~\cite{CdV85} in order to construct the subset $S_0$. Observe that
the family of subsets $(S_{\frac{1}{l}})_{l\geq 1}$ we have just defined is nonincreasing. We fix, for every $l\geq 2$,
$$\alpha_1\ \text{such that},\ \forall\lambda^2\geq\alpha_1,\ \frac{\sharp \left(S_{1}\cap\left\{j:\lambda_j^2\leq\lambda^2\right\}\right)}{N(\lambda)}\geq 1-\frac{1}{2},\ \ldots,$$
$$\alpha_l\geq\alpha_{l-1}\ \text{such that},\ \forall\lambda^2\geq\alpha_l,\ \frac{\sharp \left(S_{\frac{1}{l}}\cap\left\{j:\lambda_j^2\leq\lambda^2\right\}\right)}{N(\lambda)}\geq 1-\frac{1}{2^l}.$$ 
Then, one defines $S_0$ such that $S_0\cap\{j:\alpha_l\leq\lambda_j^2<\alpha_{l+1}\}=S_{\frac{1}{l}}\cap\{j:\alpha_l\leq\lambda_j^2<\alpha_{l+1}\}$. 
In particular, for every $l\geq 1$ and for every $\alpha_l\leq\lambda^2<\alpha_{l+1}$, one has  $S_0\cap\{j:0\leq\lambda_j^2\leq\lambda^2\}\supset S_{\frac{1}{l}}\cap\{j:0\leq\lambda_j^2\leq\lambda^2\}$ and thus, $S_0$ has density $1$.
By construction, any accumulation point of the sequence $(\mu_j(a))_{j\in S_0}$ is~$\leq A_0$.

\end{rema}

\subsection{Proof of theorem~\ref{t:GQE}}

We are now in position to prove theorem~\ref{t:GQE}. For that purpose, we interpret lemma~\ref{l:GQE} as an inequality on linear forms and then we apply Hahn-Banach Theorem.

First, we observe that
 $$\inf_{\rho\in\Lambda} L_{\rho}(a)\leq\operatorname{essinf}\ L_{\rho}(a)\leq\operatorname{esssup}\ L_{\rho}(a)\leq\sup_{\rho\in\Lambda} L_{\rho}(a).$$
Fix now $(a_k)_{k\in\mathbb{N}}$ a family of smooth functions which is dense in $\mathcal{C}^0(S^*M,\mathbb{R})$ (for the uniform topology). 
For every $k\in\mathbb{N}$, lemma~\ref{l:GQE} gives us a subset $S_{k}$ of density $1$. Without loss of generality\footnote{The intersection of two subsets of density $1$ is still a subset of density $1$.}, one can suppose that $(S_{k})_{k\in\mathbb{N}}$ is a nonincreasing sequence of subsets.
Using remark~\ref{r:cdv}, one can choose a subset $S$ of density $1$ such that 
$$\forall k\in\mathbb{N},\ \inf_{\rho\in\Lambda} L_{\rho}(a_k)\leq\liminf_{j\rightarrow+\infty,j\in S}\mu_j(a_k)\leq\limsup_{j\rightarrow+\infty,j\in S}\mu_j(a_k)\leq\sup_{\rho\in\Lambda} L_{\rho}(a_k).$$
Fix now an accumulation point $\mu$ of the sequence $(\mu_j)_{j\in S}$. By a density argument, the above inequality implies then 
$$\forall a\in\mathcal{D}(S^*M,\mathbb{R}),\ \inf_{\rho\in\Lambda} L_{\rho}(a)\leq \mu(a)\leq\sup_{\rho\in\Lambda} L_{\rho}(a).$$

As the space $\mathcal{D}(S^*M,\mathbb{R})$ is the topological dual of $\mathcal{D}'(S^*M,\mathbb{R})$ (Theorem $XIV$, Chapter~3 in~\cite{Sc66}), the previous inequality implies that, for every continuous linear form $\Phi$ on $\mathcal{D}'(S^*M,\mathbb{R})$, 
$$\inf_{\rho\in\Lambda} \Phi(L_{\rho})\leq\Phi(\mu)\leq\sup_{\rho\in\Lambda} \Phi(L_{\rho}).$$

Suppose by contradiction that $\mu$ does not belong to $\text{Cv}(\Lambda)$. By the Hahn-Banach theorem~\cite{Tr67}, there exists a continuous linear form $\Phi_0$ on $\mathcal{D}'(S^*M,\mathbb{R})$ that strictly separates the closed convex subset $\text{Cv}(\Lambda)$ from the compact convex subset $\{\mu\}$. In particular, there exists $\alpha$ in $\mathbb{R}$ such that
$$\forall \nu\in\operatorname{Cv}(L),\ \Phi_0(\nu)\leq\alpha<\Phi_0(\mu).$$
We get that $\sup_{\nu\in\operatorname{Cv}(L)}\Phi_0(\nu)\leq\alpha<\Phi_0(\mu).$ In particular, $\sup_{\rho\in \Lambda} \Phi_0(L_{\rho})\leq\alpha<\Phi_0(\mu)$ which leads to the contradiction.




\section{Proof of proposition~\ref{p:genprop}}\label{s:prop}

In this final section, we will prove proposition~\ref{p:genprop} that we used in our applications to surfaces of nonpositive curvature. 

As in the statement of the proposition, we fix an invariant subset $\Lambda'$ in $S^*M$ containing a nonempty open ball (modulo a subset of zero Liouville measure) 
and satisfying $L_{|\Lambda'}$ ergodic. 
Let $\Lambda_0\subset\Lambda'^c$ satisfying $L(\Lambda_0)=L(\Lambda'^c)$. Thanks to corollary~\ref{c:maincoro}, we know that there exists $S\subset\mathbb{N}$ of density $1$ such that any accumulation point of the sequence $(\mu_j)_{j\in S}$ is of the form
$$\alpha  \frac{L_{|\Lambda'}}{L(\Lambda')}+(1-\alpha)\nu_0,$$ 
where $\alpha\geq 0$ and where $\nu_0$ belongs to $\text{Cv}(\Lambda_0)$.\\

As $\Lambda'$ contains a nonempty open ball (modulo a set of zero Liouville measure), one can pick $\chi\geq 0$ a smooth function which is compactly supported in this open ball and which is not equal to $0$ everywhere. In particular, $\int_{S^*M}\chi dL=\int_{\Lambda'}\chi dL>0$.



\subsection{Preliminary remarks}

Recall from Birkhoff Ergodic Theorem that
$$L_{\rho}(\chi)=\lim_{T\rightarrow+\infty}\frac{1}{T}\int_0^T\chi\circ g^t(\rho)dt$$
is well defined for $L$ almost every $\rho$ in $S^*M$. From our assumptions on $\Lambda'$, one can verify that 
$$L_{\rho}(\chi)=\frac{1}{L(\Lambda')}\int_{\Lambda'}\chi dL,\ \text{a.e. on}\ \Lambda'.$$
Still thanks to Birkhoff Ergodic Theorem, we also have
$$\int_{\Lambda'}\chi dL=\int_{S^*M}\chi dL=\int_{S^*M}L_{\rho}(\chi) dL(\rho)= \int_{\Lambda'}\chi dL+\int_{\Lambda'^c}L_{\rho}(\chi)dL(\rho).$$
As $L_{\rho}(\chi)\geq 0$ a.e., we obtain that $L_{\rho}(\chi)=0$ on some $\Lambda_1\subset\Lambda'^c$ satisfying $L(\Lambda_1)=L(\Lambda'^c)$.\\ 

We pick such a subset $\Lambda_1$ in proposition~\ref{p:genprop}. Let $\Lambda_0$ be a subset of $\Lambda_1$ verifying~\eqref{e:birkhoff} and $L(\Lambda_1)=L(\Lambda_0)$. Thanks to the first part of the proposition, we know that there $S\subset\mathbb{N}$ of density $1$ such that any accumulation point of the sequence $(\mu_j)_{j\in S}$ is of the form
$$\alpha  \frac{L_{|\Lambda'}}{L(\Lambda')}+(1-\alpha)\nu_0,$$ 
where $\alpha\geq 0$ and where $\nu_0$ belongs to $\text{Cv}(\Lambda_0)$. We underline that the properties of $\Lambda_0$ and $\chi$ imply $\nu_0(\chi)=0$. Our goal is to show that there exist subsequences for which $\alpha$ can be chosen positive. For that purpose, we will study the limit of the sequence $(\mu_j(\chi))_{j\in S}$ and show that it must be positive for some subsequences of positive density.


\subsection{Proof of the proposition}

First, we observe that one can again replace $\Op$ by a nonnegative quantization procedure $\Op^+$. In particular, one gets, as $j\rightarrow+\infty$,
$$\mu_j(\chi)=\langle\psi_j,\Op^+(\chi)\psi_j\rangle+o(1),$$ 
and $\langle\psi_j,\Op^+(\chi)\psi_j\rangle\geq 0$ (as $\chi\geq 0$). Thus, without loss of generality, one can look at the accumulation points of the sequence
$$\mu_j^+(\chi):=\langle\psi_j,\Op^+(\chi)\psi_j\rangle,\ j\in S,$$
where $S\subset\mathbb{N}$ is of density $1$. Thanks to the local Weyl law~\eqref{e:weyl}, we observe that
$$\lim_{\lambda\rightarrow+\infty}\frac{1}{N(\lambda)}\sum_{j\in S:\lambda_j^2\leq\lambda^2}\mu_j^+(\chi)=\lim_{\lambda\rightarrow+\infty}\frac{1}{N(\lambda)}\sum_{j:\lambda_j^2\leq\lambda^2}\mu_j(\chi)=\int_{S^*M}\chi dL=\int_{\Lambda'}\chi dL.$$

Moreover, using lemma~\ref{l:GQE} and the preliminary remarks, one can find a subset $ S'\subset S$ of density $1$ such that any accumulation point of the sequence $(\mu_j^+(\chi))_{j\in S'}$ belongs to the interval 
$$\left[\operatorname{essinf}L_{\rho}(\chi),\operatorname{esssup}L_{\rho}(\chi)\right]=\left[0,\frac{\int_{\Lambda'}\chi dL}{L(\Lambda')}\right].$$ 

We introduce the notation
$$\alpha_j:=\frac{\mu_j^+(\chi)}{\int_{\Lambda'}\chi dL}\geq 0.$$

In order to prove our proposition, it remains to verify that, for every $0<\epsilon\leq 1$, there exists a subset $S_{\epsilon}\subset S'$ of density $\geq\frac{1-\epsilon}{L(\Lambda')^{-1}-\epsilon}$ such that any accumulation point of the subsequence $(\alpha_j)_{j\in S_{\epsilon}}$ belongs to the interval $[\epsilon,L(\Lambda')^{-1}]$. The proof is quite straightforward: we briefly explain it for the sake of completeness.


Fix now $0<\epsilon\leq 1$. Let $\eta\ll \epsilon$ be a small positive number. From the properties of the sequence $(\alpha_j)_{j\in S'}$, there exists $A>0$ (depending on $\eta$) such that, for $\lambda^2\geq A$,
$$1-\eta\leq\frac{1}{N(\lambda)}\sum_{j\in S': \lambda_j^2\leq\lambda^2}\alpha_j\ \text{and}\ (j\in S'\ \text{and}\ \lambda_j^2> A\Longrightarrow\alpha_j\leq L(\Lambda')^{-1}+\eta).$$
Thus, one gets, for $\lambda^2\geq A$,
$$1-\eta\leq\frac{1}{N(\lambda)}\sum_{j\in S': \lambda_j^2\leq A}\alpha_j+\epsilon\frac{1}{N(\lambda)}\sharp\left\{j\in S':A<\lambda_j^2\leq\lambda^2\ \text{and}\ \alpha_j<\epsilon\right\}$$
$$\hspace{5cm}+(L(\Lambda')^{-1}+\eta)\frac{1}{N(\lambda)}\sharp\left\{j\in S':A<\lambda_j^2\leq\lambda^2\ \text{and}\ \alpha_j\geq\epsilon\right\}.$$
It implies that, for $\lambda^2$ large enough (depending on $\eta$ and on $A$),
$$1-\eta\leq\eta+\epsilon+(L(\Lambda')^{-1}+\eta-\epsilon)\frac{1}{N(\lambda)}\sharp\left\{j\in S':\lambda_j^2\leq\lambda^2\ \text{and}\ \alpha_j\geq\epsilon\right\}.$$
In other words, we have shown that, for every $\eta>0$,
$$\liminf_{\lambda\rightarrow+\infty}\frac{1}{N(\lambda)}\sharp\left\{j\in S':\lambda_j^2\leq\lambda^2\ \text{and}\ \alpha_j\geq\epsilon\right\}\geq\frac{1-2\eta-\epsilon}{L(\Lambda')^{-1}+\eta-\epsilon},$$
which implies the result.



\subsection{The case of several components}\label{ss:donnayproof}

In this last paragraph, we will discuss the more general setting of remark~\ref{r:generalization}. We want to show that proposition~\ref{p:genprop} can be slightly improved in order to allow several subsets, i.e. $|I|\geq 2$ with the notations of paragraph~\ref{ss:Donnay}.

Precisely, we want to verify that $\alpha$ can be chosen $>0$ in equation~\eqref{e:donnay} for a subsequence of positive density of eigenstates (without giving precise informations on the density of the subset). 
For that purpose, we use again the fact that $\Lambda'$ contains a nonempty open ball (modulo a subset of measure $0$) and we pick $\chi\geq 0$ a nonzero smooth function compactly supported in this open ball. As above, we can find some subset $\Lambda_1\subset\Lambda'^c$ satisfying $L(\Lambda_1)=L(\Lambda'^c)$ and $L_{\rho}(\chi)=0$ on $\Lambda_1$.\\ 

Let $\Lambda_0$ be a subset of $\Lambda_1$ verifying~\eqref{e:birkhoff} and $L(\Lambda_1)=L(\Lambda_0)$. As above, we consider a subset $S\subset\mathbb{N}$ of density $1$ such that any accumulation point of $(\mu_j)_{j\in S}$ is of the form given by equation~\eqref{e:donnay}.



We can apply lemma~\ref{l:GQE} to this function: there exists $S'\subset S$ of density $1$ such that any accumulation point of the subsequence $(\mu_j(\chi))_{j\in S'}$ belongs to the interval
$$\left[0,\max_{i\in I}\frac{\int_{\Lambda_i}\chi dL}{L(\Lambda_i)}\right].$$
An application of the local Weyl law~\eqref{e:weyl} tells us that
$$\lim_{\lambda\rightarrow+\infty}\frac{1}{N(\lambda)}\sum_{j\in S':\lambda_j^2\leq\lambda^2}\mu_j(\chi)=\lim_{\lambda\rightarrow+\infty}\frac{1}{N(\lambda)}\sum_{j:\lambda_j^2\leq\lambda^2}\mu_j(\chi)=\int_{S^*M}\chi dL>0.$$
As $\mu_j(\chi)$ is ``almost positive'', we can apply the argument of the previous paragraph. In particular, for every $\epsilon>0$ small enough, we obtain the existence of a subset $S_{\epsilon}\subset S'$ of positive density such that any accumulation point of the sequence $(\mu_j(\chi))_{j\in S_{\epsilon}}$ belongs to the interval
$$\left[\epsilon,\max_{i\in I}\frac{\int_{\Lambda_i}\chi dL}{L(\Lambda_i)}\right].$$

As $S_{\epsilon}\subset S$, one knows that any accumulation point of $(\mu_j)_{j\in S_{\epsilon}}$ is of the form
$$\mu=\alpha\sum_{i\in I}t_i\frac{L_{|\Lambda_i}}{L(\Lambda_i)}+(1-\alpha)\nu_0,$$
where $0\leq \alpha\leq 1$, $0\leq t_i\leq 1$, $\sum_{i\in I}t_i=1$ and $\nu_0$ belongs to $\text{Cv}(\Lambda_0)$. Finally, as $\mu(\chi)\geq\epsilon$ and $\nu_0(\chi)=0$ (as $L_{\rho}(\chi)=0$ on $\Lambda_0\subset\Lambda_1$), we find that $\alpha$ must be positive.


\section*{Acknowledgements} This work has been partially supported by the Labex CEMPI (ANR-11-LABX-0007-01) and by the grant ANR-09-JCJC-0099-01 of the Agence Nationale de la Recherche. We warmly thank Stephan De Bi\`evre for discussions and comments related to this work, and the referee for precious advices in order to improve the presentation of this article.

\end{document}